\documentclass[10pt,a4paper]{amsart}
\RequirePackage[OT1]{fontenc}
\RequirePackage{amsthm,amsmath,amssymb,mathrsfs, enumerate}
\RequirePackage[colorlinks,citecolor=blue,urlcolor=blue]{hyperref}
\RequirePackage{hyphenat}
\usepackage{tikz}

\newcommand{\be}{\begin{equation}}
\newcommand{\ee}{\end{equation}}
\newcommand{\la}{\langle}
\newcommand{\ra}{\rangle}

\def\RR{\mathbb{R}}
\def\NN{\mathbb{N}}
\def\ZZ{\mathbb{Z}}
\def\i{{\bf i}}\def\j{{\bf j}}\def\x{{\bf x}}
\def\k{{\bf k}}
\def\y{{\bf y}}

\title[Realizability problem on the $d-$dimensional lattice]
{Translation invariant realizability problem \\ on the $d-$dimensional lattice: \\ an explicit construction}
\author[E. Caglioti, M. Infusino, T. Kuna]{E. Caglioti$^{1}$, M. Infusino$^{2}$, T. Kuna$^{3}$}\address{$^{1}$ \emph{Emanuele Caglioti}\newline \indent
Dipartimento di Matematica,\newline \indent
Sapienza Universit\`a di Roma,\newline \indent 
00185, Rome, Italy.}
\email{caglioti@mat.uniroma1.it }
\address{$^{2}$ \emph{Maria Infusino}\newline \indent
Fachbereich Mathematik und Statistik,\newline \indent
Universit\"at Konstanz,\newline \indent
78457 Konstanz, Germany.}
\email{maria.infusino@uni-konstanz.de}
\address{$^{3}$ \emph{Tobias Kuna} \newline \indent
Department of Mathematics and Statistics,\newline \indent
University of Reading,\newline \indent
RG66AX, Reading, UK.}
\email{t.kuna@reading.ac.uk}
\date{}
\begin{document}
\begin{abstract}
We consider a particular instance of the truncated realizability problem on the $d-$dimensional lattice. Namely, given two functions $\rho_1({\bf i})$ and $\rho_2({\bf i},{\bf j})$ non-negative and symmetric on $\mathbb{Z}^d$, we ask whether they are the first two correlation functions of a translation invariant point process. We provide an explicit construction of such a realizing process for any $d\geq 2$ when the radial distribution has a specific form. We also derive from this construction a lower bound for the maximal realizable density and compare it with the already known lower bounds.
\end{abstract} 

\maketitle

\section{Introduction}
Let $d$ be a positive integer. Given a point process $P=\{P_\i\}_{\i\in\ZZ^d}$ on the $d-$dimen\-sional lattice, $P_\i\in\{0,1\}$, whose distribution is described by the probability measure~$\mu$, we define the first and second order \emph{correlation function} as follows
$$\left\{\begin{array}{ll}
\rho_1(\i):=\la P_\i\ra\\
\rho_2(\i, \j):=\la P_\i P_\j \ra-\rho_1(\i)\delta(\i-\j)\\
\end{array}\right.,
$$
where $\i,\j\in\ZZ^d$, $\delta$ is the Dirac delta function and $\la\cdot\ra$ denotes the expectation w.r.t.~$\mu$.

The \emph{truncated realizability problem} addresses the inverse question. Namely, given two functions $\rho_1(\i)$ and $\rho_2(\i,\j)$ non-negative and symmetric for all $\i,\j\in\ZZ^d$, does there exist a point process $P$ for which these are the correspondent first and second order correlation functions? Clearly, the truncated realizability problem can be posed for any finite sequence of non-negative and symmetric functions $(\rho_k(\i_1,\i_2,\ldots,\i_k))_{k=1}^n$. When the question is asked for a given infinite sequence $(\rho_k)_{k\in\NN}$, then the problem is addressed as full realizability problem (see e.g. \cite{Le75a, Le75b} for a systematic study of the full realizability problem for point processes and \cite{IK, IKR14} for recent developments).

In the following we will consider the important special case of \emph{translation invariant point processes}, which actually contains all the essential difficulties of the problem. In this case the realizability problem asks if, for given $\rho\in\RR^+$ and $g:\ZZ^d\to\RR^+$ symmetric, there exists a translation invariant point process such that its first two correlation functions are given by
\be\label{transl-inv}
\left\{\begin{array}{ll}
\rho_1(\i)=\rho\\
\rho_2(\i, \j)=\rho^2g(\i-\j)\\
\end{array}\right..
\ee
If such a process exists, then it is said to be \emph{realizing} and the pair $(\rho, g)$ is called \emph{realizable} on $\ZZ^d$. Note that writing the second order correlation in this form is not an additional restriction beyond the assumption of translation invariance. The function $g$ is known in classic fluid theory as \emph{radial distribution},~\cite{HaMcDo87}. 

The truncated realizability problem is in fact a longstanding problem in the classical theory of fluids (see e.g. \cite{GaPerc64, Perc64, PercUnpub}), but it has been investigated in many other contexts such as stochastic geometry~\cite{Mol05},  spatial statistics~\cite{CoMon12, Stoy00}, spatial ecology~\cite{MDL04} and neural spike trains~\cite{BKM04, JM04}, just to name a few. In particular, Stillinger, Torquato et al. developed fascinating applications in the study of heterogeneous materials and mesoscopic structures based on the solvability of the truncated realizability problem (see e.g.\! \cite{CraToSt03, ToSt04, ToSt05, To02, ToSt03, ToSt06}). A structural investigation of this problem was recently started in \cite{KuLeSp11}, where the authors identify the realizability problem as a particular instance of the infinite-dimensional truncated moment problem (see \cite{CaKuLeSp06, Kor05, Kor07, LeKuSp07, MolLach15} for further recent developments about the truncated realizability problem for point processes). As far as we know, the only earlier reference about the truncated infinite-dimensional moment problem is~\cite{Us74}.

In this paper, we will show how to explicitly construct a point process on the $d-$dimensional lattice with $d\geq 2$ such that, for given $\alpha\geq 0$, \eqref{transl-inv} holds for certain values of $\rho$ and for $g=g^{(\alpha)}$ defined as follows: 
\be\label{galpha}
g^{(\alpha)}(\x):=\left\{\begin{array}{ll}
0 &\text{if } \,\,\x\,=0\\
\alpha &\text{if } |\x|=1\\
1 &\text{if } |\x|>1\\
\end{array}\right..
\ee
Explicit constructions of point processes realizing this lattice problem in the case $d=1$ were provided in \cite[Appendix~1]{LeKuSp07}.
Such a problem has been extensively studied for the case $\alpha=0$ by Stillinger and Torquato in \cite{ToSt04} (see also \cite{CraToSt03, ToSt03}). The function $g^{(0)} $ describes a model with on-site and nearest neighbour exclusion and with no correlation for pairs of sites separated by two or more lattice spacings.

From \cite[Section~1]{LeKuSp07}, we know that for a fixed $\alpha$ the set of realizable densities $\rho$ is an interval $[0,\bar{\rho}_\alpha(d)]$ with $0<\bar{\rho}_\alpha(d)\leq 1$. Moreover, in \cite{LeKuSp07} the authors discuss:
\begin{enumerate}[(i)]
\item certain general methods which, when applied to \eqref{galpha}, yield lower bounds for $\bar\rho_\alpha(d)$ in any dimension $d$.
\item concrete upper and lower bounds for $\bar\rho_\alpha(1)$. In particular, the lower bounds improve those obtained from the general methods~(i).
\end{enumerate}
Our $d-$dimensional construction combined with the one-dimensional lower bounds (ii) provides a lower bound for $\bar{\rho}_\alpha(d)$ for any $d\geq 2$ and any $\alpha\geq 0$.  We will briefly compare this with the lower bound obtained from the general methods~(i). We also follow techniques from \cite{LeKuSp07} to get an upper bound for $\bar{\rho}_\alpha(d)$.

\section{An explicit realizing translation invariant point process on~$\ZZ^d$}\label{ExplConstr}
In the following, we explicitly construct a point process $P=\{P_\i\}_{\i\in\ZZ^d}$ on the $d-$dimensional lattice with $d\geq 2$ such that, for given $\alpha\geq 0$, the following hold for certain values of $\rho$ (depending on $\alpha$ and on $d$):
\be\label{R1}
\la P_\i\ra=\rho
\ee and
\be\label{R2}
\la P_\i P_\j\ra=\left\{\begin{array}{ll}
\rho &\text{if } \,\i=\j\\
\alpha\rho^2 &\text{if } |\i-\j|=1\\
\rho^2 &\text{if } |\i-\j|>1\\
\end{array}
\right.,
 \ee
 that is, the radial distribution is given by \eqref{galpha}.
\subsection{Construction in dimension 2}\label{ExplConstrDim2}\ \vspace{-0.3cm} \\

In order to build such a process on $\ZZ^2$ we start from a realizing one-dimensional process achieving density $\gamma$. Namely, given $\alpha\geq 0$, we consider a point process $\{A_i\}_{i\in\ZZ}$, $A_i\in\{0,1\}$, defined on the one-dimensional lattice and such that for some $\gamma>0$ we have
$$\la A_i\ra=\gamma,$$ and
$$
\la A_iA_j\ra=\left\{\begin{array}{ll}
\gamma &\text{if } \,i=j\\
\alpha\gamma^2 &\text{if } |i-j|=1\\
\gamma^2 &\text{if } |i-j|>1\\
\end{array}
\right..
$$
We denote a process of this kind by $BP\gamma$ that stays for \emph{basic process with density~$\gamma$}. As pointed out in the introduction, there exists a good number of constructions of realizing processes in the one-dimensional case, see e.g.\! \cite[Appendix~1]{LeKuSp07}. The results in the one-dimensional case relevant to our investigation (in particular the range where $\gamma$ can vary) are recalled in Section~3, which is devoted to the discussion of the maximal realizable density in any dimension.

Let us define  two processes $B^{(1)}=\left\{B_{i_1,i_2}^{(1)}\right\}_{(i_1,i_2)\in\ZZ^2}$ and $B^{(2)}=\left\{B_{i_1,i_2}^{(2)}\right\}_{(i_1,i_2)\in\ZZ^2}$ on $\ZZ^2$ as follows. For a fixed $i_1\in\ZZ$, the process $\left\{B_{i_1,i_2}^{(1)}\right\}_{i_2\in\ZZ}$ is a $BP\gamma$ in~$i_2$. For any $i_1,j_1\in\ZZ$ with $i_1\neq j_1$,  the processes $\left\{B_{i_1,i_2}^{(1)}\right\}_{i_2\in\ZZ}$ and $\left\{B_{j_1,j_2}^{(1)}\right\}_{j_2\in\ZZ}$ are independent. In particular, we have
\be\label{tre}
\la B_{i_1,i_2}^{(1)}B_{j_1,j_2}^{(1)}\ra=\left\{\begin{array}{ll}
\gamma^2 &\text{if } i_1\neq j_1\\
\gamma &\text{if } i_1= j_1 \text{ and } i_2=j_2\\
\alpha\gamma^2 &\text{if } i_1= j_1 \text{ and } |i_2-j_2|=1\\
\gamma^2 &\text{if } i_1= j_1 \text{ and } |i_2-j_2|>1\\
\end{array}
\right..
 \ee
In other words, the process $B^{(1)}$ can be seen as a sequence of vertical $BP\gamma$'s independent one from each other (see Figure 1 for an example). 

 \begin{figure}[h]
 \centering
 \begin{tikzpicture}
\draw[step=1cm,gray,very thin] (0,0) grid (5,5);
\draw[ultra thick,-] (0,0) -- (0,5);
\draw[ultra thick,-] (1,0) -- (1,5);
\draw[ultra thick,-] (2,0) -- (2,5);
\draw[ultra thick,-] (3,0) -- (3,5);
\draw[ultra thick,-] (4,0) -- (4,5);
\draw[ultra thick,-] (5,0) -- (5,5);
\node[circle,draw, thick, fill=black!20] at (0,1) {};
\node[circle,draw, thick, fill=black!20] at (0,3) {};
\node[circle,draw, thick, fill=black!20] at (1,0) {};
\node[circle,draw, thick, fill=black!20]  at (1,2) {};
\node[circle,draw, thick, fill=black!20]  at (1,4) {};
\node[circle,draw, thick, fill=black!20]  at (2,1) {};
\node[circle,draw, thick, fill=black!20]  at (2,3) {};
\node[circle,draw, thick, fill=black!20]  at (2,5) {};
\node[circle,draw, thick, fill=black!20]  at (3,0) {};
\node[circle,draw, thick, fill=black!20]  at (3,3) {};
\node[circle,draw, thick, fill=black!20]  at (3,5) {};
\node[circle,draw, thick, fill=black!20]  at (4,1) {};
\node[circle,draw, thick, fill=black!20]  at (4,3) {};
\node[circle,draw, thick, fill=black!20]  at (5,1) {};
\node[circle,draw, thick, fill=black!20]  at (5,5) {};
\end{tikzpicture}
\caption{Example of process $B^{(1)}$ (with $\alpha=0$)}
\end{figure}
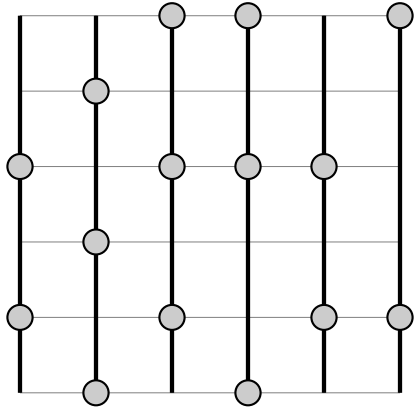
Similarly, the process $B^{(2)}$ is defined as a sequence of horizontal $BP\gamma$'s independent one from each other (see Figure 2 for an example), i.e.
\be\label{quattro}
\la B_{i_1,i_2}^{(2)}B_{j_1,j_2}^{(2)}\ra=\left\{\begin{array}{ll}
\gamma^2 &\text{if } i_2\neq j_2\\
\gamma &\text{if } i_2= j_2 \text{ and } i_1=j_1\\
\alpha\gamma^2 & \text{if } i_2= j_2 \text{ and } |i_1-j_1|=1\\
\gamma^2 &\text{if } i_2= j_2 \text{ and } |i_1-j_1|>1\\
\end{array}
\right..
 \ee
 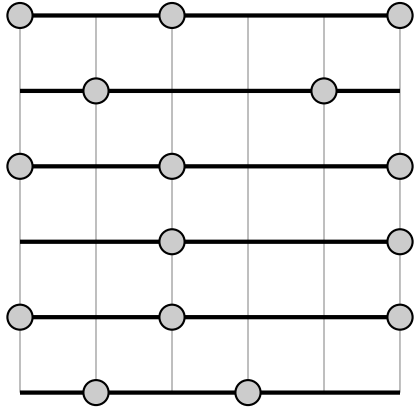
\begin{figure}[h]
  \centering
 \begin{tikzpicture}
\draw[step=1cm,gray,very thin] (0,0) grid (5,5);
\draw[ultra thick,-] (0,0) -- (5,0);
\draw[ultra thick,-] (0,1) -- (5,1);
\draw[ultra thick,-] (0,2) -- (5,2);
\draw[ultra thick,-] (0,3) -- (5,3);
\draw[ultra thick,-] (0,4) -- (5,4);
\draw[ultra thick,-] (0,5) -- (5,5);
\node[circle,draw, thick, fill=black!20]  at (1,0) {};
\node[circle,draw, thick, fill=black!20]  at (3,0) {};
\node[circle,draw, thick, fill=black!20]  at (0,1) {};
\node[circle,draw, thick, fill=black!20]  at (0,3) {};
\node[circle,draw, thick, fill=black!20]  at (2,1) {};
\node[circle,draw, thick, fill=black!20]  at (5,1) {};
\node[circle,draw, thick, fill=black!20]  at (2,2) {};
\node[circle,draw, thick, fill=black!20]  at (5,2) {};
\node[circle,draw, thick, fill=black!20]  at (2,3) {};
\node[circle,draw, thick, fill=black!20]  at (5,3) {};
\node[circle,draw, thick, fill=black!20]  at (1,4) {};
\node[circle,draw, thick, fill=black!20]  at (4,4) {};
\node[circle,draw, thick, fill=black!20]  at (0,5) {};
\node[circle,draw, thick, fill=black!20]  at (2,5) {};
\node[circle,draw, thick, fill=black!20]  at (5,5) {};
\end{tikzpicture}
\caption{Example of process $B^{(2)}$ (with $\alpha=0$)}
\end{figure}

Let us define now the process $P=\{P_\i\}_{\i\in\ZZ^2}$ as
$$P_{i_1,i_2}:=B_{i_1,i_2}^{(1)}B_{i_1,i_2}^{(2)},$$
(see Figure~3 for $P$ constructed from the basic processes in Figures~1 and 2). Since the processes $B^{(1)}$ and $B^{(2)}$ are independent, we get 
$$\la P_{i_1,i_2}\ra=\la B_{i_1,i_2}^{(1)}B_{i_1,i_2}^{(2)}\ra=\la B_{i_1,i_2}^{(1)}\ra\la B_{i_1,i_2}^{(2)}\ra=\gamma\cdot\gamma=\gamma^2.$$
Hence, \eqref{R1} holds for $\rho=\gamma^2$.
From the independence of $B^{(1)}$ and $B^{(2)}$, we also get 
\be\label{star}
\la P_{i_1,i_2}P_{j_1,j_2}\ra=\la B_{i_1,i_2}^{(1)}B_{i_1,i_2}^{(2)}B_{j_1,j_2}^{(1)}B_{j_1,j_2}^{(2)}\ra=\la B_{i_1,i_2}^{(1)}B_{j_1,j_2}^{(1)}\ra\la B_{i_1,i_2}^{(2)}B_{j_1,j_2}^{(2)}\ra.
\ee
We can easily check, by using \eqref{tre} and \eqref{quattro} in \eqref{star}, that \eqref{R2} holds for $\rho=\gamma^2$. In fact, we need to consider only the following four cases, because all the others are equivalent to these ones by symmetry.
\begin{enumerate}[a)]
\item If $i_1=j_1$ and $i_2=j_2$ then  $$\la P_{i_1,i_2}P_{j_1,j_2}\ra=\gamma\cdot\gamma=\gamma^2.$$
\item If $i_1=j_1$ and $i_2\neq j_2$ then $$\la P_{i_1,i_2}P_{j_1,j_2}\ra=\la B_{i_1,i_2}^{(1)}B_{j_1,j_2}^{(1)}\ra\gamma^2.$$
Therefore:
\begin{itemize}
\item if $i_2=j_2+1$ then $\la B_{i_1,i_2}^{(1)}B_{j_1,j_2}^{(1)}\ra=\alpha\gamma^2$ and so $\la P_{i_1,i_2}P_{j_1,j_2}\ra=\alpha\gamma^4$\vspace{0.2cm}
 \item if $|i_2-j_2|>1$ then $\la B_{i_1,i_2}^{(1)}B_{j_1,j_2}^{(1)}\ra=\gamma^2$ and so $\la P_{i_1,i_2}P_{j_1,j_2}\ra=\gamma^4$.
 \end{itemize}\vspace{0.2cm}
 \item If $i_1=j_1+1$ and $|i_2-j_2|>1$ then 
 $$\la P_{i_1,i_2}P_{j_1,j_2}\ra=\la B_{i_1,i_2}^{(1)}B_{j_1,j_2}^{(1)}\ra\la B_{i_1,i_2}^{(2)}B_{j_1,j_2}^{(2)}\ra=\gamma^4.$$
 \item If $|i_1-j_1|>1$ and $i_2\neq j_2$ then 
 $$\la P_{i_1,i_2}P_{j_1,j_2}\ra=\la B_{i_1,i_2}^{(1)}B_{j_1,j_2}^{(1)}\ra\la B_{i_1,i_2}^{(2)}B_{j_1,j_2}^{(2)}\ra=\gamma^4.$$
\end{enumerate}
 \begin{figure}[h]
  \centering
 \begin{tikzpicture}
\draw[step=1cm,gray,very thin] (0,0) grid (5,5);
\node[circle,draw, thick, fill=black!20] at (0,1) {};
\node[circle,draw, thick, fill=black!20] at (0,3) {};
\node[circle,draw, thick, fill=black!20]  at (1,4) {};
\node[circle,draw, thick, fill=black!20]  at (2,1) {};
\node[circle,draw, thick, fill=black!20]  at (2,3) {};
\node[circle,draw, thick, fill=black!20]  at (2,5) {};
\node[circle,draw, thick, fill=black!20]  at (3,0) {};
\node[circle,draw, thick, fill=black!20]  at (5,1) {};
\node[circle,draw, thick, fill=black!20]  at (5,5) {};
\node[circle,draw, thick, fill=black!20]  at (1,0) {};
\end{tikzpicture}
\caption{Process $P$ constructed from the processes in Figure~1 and Figure~2 }
\end{figure}
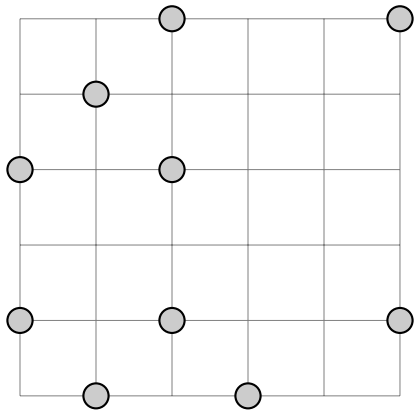

\subsection{Construction in higher dimension}\label{ExplConstrDimd}\ \vspace{-0.3cm}\\

The construction presented in the previous subsection easily generalizes to any dimension $d>2$ by defining $P_{i_1,\ldots,i_d}:=B^{(1)}_{i_1,\ldots,i_d}\cdots B^{(d)}_{i_1,\ldots,i_d}$
where, for any fixed $i_2,\ldots,i_d\!\in\!\ZZ$, $\{B^{(1)}_{i_1,\ldots,i_d}\}_{i_1\in\ZZ}$ is a $BP\gamma$ in the variable $i_1$ with density~$\gamma$ and similarly for the other variables.  Therefore, the point process $P$ on $\ZZ^d$ defined as above satisfies \eqref{R1} and \eqref{R2} for $\rho=\gamma^d$.

\section{Bounds for the maximal realizable density}
In this section, we will discuss the problem of estimating the maximal realizable density~$\bar{\rho}_\alpha(d)$. In particular, we will show a general upper bound for any $d\geq 1$ using the technique introduced in \cite{LeKuSp07} for $d=1$. As for the lower bound, we will recall the results in \cite[Appendix~1]{LeKuSp07} for the one-dimensional case and combine them with the explicit construction proposed in Section~\ref{ExplConstr} to produce a lower bound for~$\bar{\rho}_\alpha(d)$ for any $d\geq 2$. We will compare this with the lower bound obtained by applying the general methods of \cite{LeKuSp07} to the case when the radial distribution is given by \eqref{galpha}.
\subsection{Upper bounds for $\bar{\rho}_\alpha(d)$}\ \vspace{-0.3cm}\\

For $d\geq 1$ and $\alpha\geq 0$, the covariance matrix associated to a given pair $(\rho, g^{(\alpha)})$ realizable on $\ZZ^d$ must be positive semidefinite. This is equivalent to the non-negativity of the corresponding infinite volume structure function $\hat{S}$ on $\RR^d$ (for more details see e.g. \cite[Section 2]{LeKuSp07}):
$$\hat{S}(\k):=\rho+\rho^2\sum_{\x\in\ZZ^d}e^{i\k\cdot\x}[g^{(\alpha)}(\x)-1]\geq 0,\,\forall\,\k\in\RR^d.$$
This leads to an explicit upper bound for the maximal realizable density $\bar{\rho}_\alpha(d)$. 
In fact, it is easy to see that for any $\k:=(k_1,\ldots, k_d)\in\RR^d$ we get
\begin{align*}
\hat{S}(\k)&=\rho-\rho^2+\rho^2\sum_{\x\in\ZZ^d, |\x|=1}e^{i\k\cdot\x}(\alpha-1)\\
&=\rho-\rho^2+\rho^2(\alpha-1)\sum_{j=1}^d (e^{ik_j}+e^{-ik_j})\\
&=\rho\left[1-\rho\left(1-2(\alpha-1)\sum_{j=1}^d \cos(k_j)\right)\right].
\end{align*}
Then, using the non-negativity of $\hat{S}$ on $\RR^d$, we get that 
$$\rho\leq \frac 1{f_\alpha(k_1,\ldots, k_d)},\quad\forall\,(k_1,\ldots, k_d)\in\RR^d,$$
where
$f_\alpha(k_1,\ldots, k_d):=1-2(\alpha-1)\sum_{j=1}^d \cos(k_j).$
The best upper bound is then obtained for the points of $\RR^d$ where $f_\alpha$ attains the maximum. Hence, we have that
$$\bar{\rho}_\alpha(d)\leq \frac 1{\max\limits_{\k\in\RR^d}f_\alpha(\k)}=:R_F(\alpha, d).$$
By computing the maximum of $f_\alpha$ over $\RR^d$, we get our upper bound
\be\label{UB}
R_F(\alpha, d)=\frac{1}{1+2d|1-\alpha|}.
\ee

As mentioned above, this technique was employed in \cite[Appendix 1]{LeKuSp07} to get $R_F(\alpha,1)$. Furthermore, the authors provided another upper bound $R_Y(\alpha,1)$ in the one-dimensional case by using the Yamada condition (see \cite{Yam61}). Note that 
$$
\left\{\begin{array}{ll}
R_Y(\alpha,1)=R_F(\alpha,1), & \text{if }\alpha=\frac 12\text{ or }\alpha=\frac{k\pm1}{2k}, k\in\NN \text{ or }\alpha\geq 1\\
R_Y(\alpha,1)<R_F(\alpha,1) & \text{otherwise}.
\end{array}\right. $$

\subsection{Lower bounds for $\bar{\rho}_\alpha(d)$}\ \vspace{-0.3cm} \\

Applying \cite[Theorem~3.2]{LeKuSp07} for $g\equiv g^{(\alpha)}$ when $0\leq \alpha<1$ and \cite[Theorem~5.1]{LeKuSp07} for $G_2(\x,\y)=g^{(\alpha)}(\y-\x)$ for all $\x,\y\in\ZZ^d$ such that $\x\neq\y$ when $\alpha\geq 1$, we get that \be\label{d-dim-ArmLB}
\bar{\rho}_\alpha(d)\geq r_A(\alpha,d):=\left\{\begin{array}{ll}
\frac{1}{e(2d+1-2d\alpha)}, & \text{if }0\leq\alpha< 1,\\
\frac{1}{\alpha^{2d}}, & \text{if }\alpha\geq 1.
\end{array}\right. 
\ee
For $d=1$, this lower bound has been improved in \cite[Appendix 1]{LeKuSp07} by explicitly constructing a translation invariant realizing process at some value of $\rho$ and $\alpha$. Let us summarize in one formula the lower bounds coming from the two main constructions considered in \cite[Appendix 1]{LeKuSp07}: 
\be\label{1-dimExplLB}
\bar{\rho}_\alpha(1)\geq \left\{\begin{array}{ll}
\frac{1}{(1+\sqrt{1-\alpha})^2}, & \text{if } 0\leq\alpha< \frac 12,\\
\frac{1}{1+\sqrt{2-2\alpha}}, & \text{if }\frac 12\leq\alpha\leq 1,\\
\frac{1}{2\alpha-1}, & \text{if }\alpha\geq 1.\\
\end{array}\right.
\ee
In \cite{CaKuLeSp06} a further explicit construction is provided for the case $\alpha=0$, which slightly improves this lower bound to $\bar{\rho}_0(1)>0.265$. In the same work also the upper bound is improved to $\bar{\rho}_0(1)<(326-\sqrt{3115})/822\approx 0.3287$. However, it remains an open problem to reduce the gap between lower and upper bounds for $\bar{\rho}_\alpha(1)$ for any $\alpha\geq 0$.

Exactly as in the one-dimensional case, also for $d\geq 2$, one can try to obtain better lower bounds than \eqref{d-dim-ArmLB} by using explicit constructions. In the following, we will use the construction we proposed in Section~\ref{ExplConstr} combined with the one-dimensional lower bound \eqref{1-dimExplLB} to compute a new lower bound for $\bar{\rho}_\alpha(d)$, which we will briefly compare with \eqref{d-dim-ArmLB}.

If we apply the construction given for $d\geq 2$ in Section~\ref{ExplConstr} starting with a basic process with density $\bar{\rho}_\alpha(1)$, then we get a point process on $\ZZ^d$ which realizes the pair $\left((\bar{\rho}_\alpha(1))^d, g^{(\alpha)}\right)$ for any $\alpha\geq 0$. This explicit construction guarantees that for any $\alpha\geq 0$, $$\bar{\rho}_\alpha(d)\geq(\bar{\rho}_\alpha(1))^d.$$
Using the lower bounds \eqref{1-dimExplLB} in the latter inequality, we directly have the following
\be\label{NewLowerBound}
\bar{\rho}_\alpha(d)\geq r_C(\alpha,d):= \left\{\begin{array}{ll}
\frac{1}{(1+\sqrt{1-\alpha})^{2d}}, & \text{if } 0\leq\alpha< \frac 12,\\
\frac{1}{(1+\sqrt{2-2\alpha})^d}, & \text{if }\frac 12\leq\alpha\leq 1,\\
\frac{1}{(2\alpha-1)^d}, & \text{if }\alpha\geq 1.\\
\end{array}\right.
\ee
Note that:
\begin{itemize}
\item if $0\leq\alpha< \frac 12$ then $r_C(\alpha,d)\leq r_A(\alpha,d)$
\item if $\alpha\geq 1$ then $r_A(\alpha,d)\leq r_C(\alpha,d)$
\item if $\frac 12\leq\alpha\leq 1$ then the relation between the two bounds depends on the dimension~$d$. Actually, for each $d\geq 2$ there exists $\alpha_C(d)\in\left[\frac 12, 1\right]$ such that $r_A(\alpha,d)\leq r_C(\alpha,d)$ for any $\alpha_C(d)\leq\alpha\leq 1$.
\end{itemize}
The comparison between the lower bounds $r_C(\alpha,d)$ and $r_A(\alpha,d)$ is illustrated in Figure~4 for $d=2,\ldots,6$ and for $0\leq\alpha<1$.

\begin{figure}[h]
\includegraphics[width=14cm]{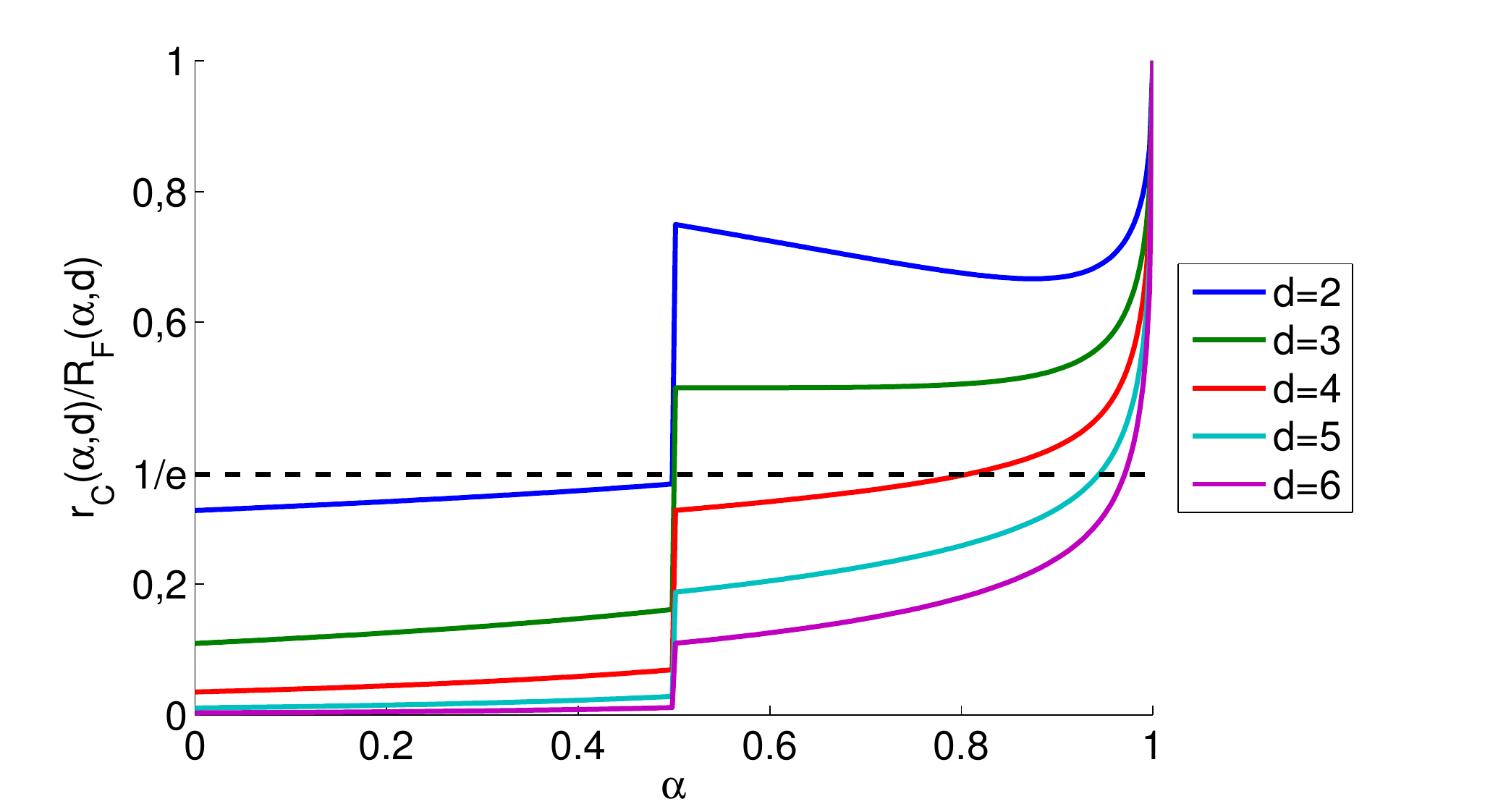}
\caption{Comparison between the lower bounds $r_C(\alpha,d)$ and $r_A(\alpha,d)$ plotted relatively to the upper bound $R_F(\alpha,d)$ as functions of $\alpha$ with $0\leq\alpha<1$. The coloured lines correspond to $\frac{r_C(\alpha,d)}{R_F(\alpha,d)}$ for $d=2,\ldots, 6$ and the dotted line to $\frac{r_A(\alpha,d)}{R_F(\alpha,d)}$ for any $d$. For the definitions of $R_F(\alpha,d)$, $r_A(\alpha,d)$, $r_C(\alpha,d)$ see \eqref{UB}, \eqref{d-dim-ArmLB}, \eqref{NewLowerBound}, respectively.
}
\end{figure}




\bibliographystyle{plain}

\section*{Acknowledgments}
The authors would like to thank Joel Lebowitz and Eugene Speer for helpful comments on the manuscript and inspiring discussions. We are very grateful to Davide Faranda for improving the graphics.

The work of Maria Infusino and Tobias Kuna was partially supported by the EPSRC Research Grant EP/H022767/1. The work of Maria Infusino was also supported by a Marie Curie fellowship of the Istituto Nazionale di Alta Matematica (Grant PCOFUND-GA-2009-245492). 	\\

\end{document}